\documentclass[11pt,a4paper,leqno]{article}
\usepackage{amsmath}
\usepackage{amssymb}
\advance\topmargin by -2.2cm
\newcommand{\bbr}{I\!\!R}

\newcommand{\bbn}{I\!\!N}

\newcommand{\lala}{\lambda\!\!\lambda}

\newcommand{\calb}{{\cal B}}
\newcommand{\calc}{{\cal C}}
\newcommand{\cald}{{\cal D}}
\newcommand{\cale}{{\cal E}}

\newcommand{\call}{{\cal L}}

\newcommand{\caln}{{\cal N}}

\newcommand{\barr}{\begin{array}}
\newcommand{\earr}{\end{array}}
\newcommand{\beqq}{\begin{equation}}
\newcommand{\eeqq}{\end{equation}}
\newcommand{\beao}{\begin{eqnarray*}}
\newcommand{\eeao}{\end{eqnarray*}\noindent}
\newcommand{\beam}{\begin{eqnarray}}
\newcommand{\eeam}{\end{eqnarray}\noindent}

\newcommand{\halmos}{\quad\hfill\mbox{$\Box$}}

\newcommand{\la}{\lambda}
\newcommand{\La}{\Lambda}

\newcommand{\si}{\sigma}
\newcommand{\al}{\alpha}

\newcommand{\vep}{\varepsilon}

\newcommand{\wt}{\widetilde}

\newcommand{\lra}{\longrightarrow}

\begin{document}

{\huge On ergodicity properties for time inhomogeneous Markov processes with $T$-periodic semigroup }  

{\bf Reinhard H\"opfner}, Universit\"at Mainz,\\
{\bf Eva L\"ocherbach}, Universit\'e de Cergy-Pontoise\\

\vskip0.3cm
{\small 
{\bf Abstract: } We consider a time inhomogeneous strong Markov process $(\xi_t)_{t\ge 0}$ taking values in a Polish state space whose semigroup has a $T$-periodic structure. After reviewing some conditions which imply ergodicity of the grid chain $(\xi_{kT})_{k\in\bbn_0}$, and thus ergodicity of the $T$-segment chain $( (\xi_{kT+s})_{0\le s\le T} )_{k\in\bbn_0}$, we formulate a new  condition for $d-$dimensional diffusions. It can be easily verified in terms of drift and diffusion coefficient of the process, and allows to deal both with unbounded coefficients and possibly degenerate diffusion term. 


\vskip0.2cm
{\bf MSC:  \quad    60 J 60, 60 J 25.}
}

\vskip1.5cm
Consider a time inhomogeneous strong Markov process $(\xi_t)_{t\ge 0}$ with c\`adl\`ag paths, taking values in a Polish state space $(E,\cale)$, with semigroup $P_{s,t}(\cdot,\cdot)$, $0\le s<t<\infty$, and a distinguished $\si$-finite measure $\Lambda$ on $(E,\cale)$ such that $\cale{\otimes}\cale$--measurable densities 
$$  
P_{s,t}(x,dy) \;=\; p_{s,t}(x,y)\, \Lambda(dy) \quad,\quad 0\le s<t<\infty \;,\; x,y\in E  
$$  
exist with respect to $\Lambda$. Assume  $T$-periodicity of the semigroup  
\beqq\label{semigroupTperiodic}
p_{s,t}(x,y) \;=\; p_{kT+s,kT+t}(x,y)  \quad\mbox{for all $k$, all $x,y$, all $0\le s<t<\infty ,$}  
\eeqq
and write $i_T(t)$ for $t \,\mbox{modulo}\, T$. Diffusion processes with $T$-periodic semigroup appear e.g.\ when a particle is moving in a periodically changing double-well potential (\cite{HI-05}), in problems coming from signal transmission (e.g.\ \cite{IH-81}, \cite{G-88}, \cite{HK-11}) which do have some interesting statistical properties, and in other contexts. Consider the Markov chain of $T$-segments in the path of $\xi$ 
$$
\mathbb{X} = (\, \mathbb{X}_k \,)_{k\in\bbn_0}  \quad,\quad  \mathbb{X}_k = (\xi_{kT+s})_{0\le s\le T}
$$
taking values in $(\mathbb{D}_T,\cald_T)$, the Skorohod space of c\`adl\`ag functions $[0,T]\to E$ equipped with its Borel $\si$-field $\cald_T$ (again a Polish space), and the $T$-grid chain 
$$
X = ( X_k )_{k\in\bbn_0} \quad,\quad X_k := \xi_{kT}  \;. 
$$
The assumption (\ref{semigroupTperiodic}) on $T$-periodicity of the semigroup of $(\xi)_{t\ge 0}$ implies that both $\mathbb{X}$ and $X$ are time homogeneous Markov chains. One can show that positive Harris recurrence of the $T$-grid chain $X$ implies positive Harris recurrence of the $T$-segment chain $\mathbb{X}$ (cf.\ H\"opfner and Kutoyants \cite{HK-10a}, theorem 2.1, same proof on $(\mathbb{D}_T,\cald_T)$ as on the path space $(\mathbb{C}_T,\calc_T)$ of continuous functions $[0,T]\to E$); in fact both properties are equivalent. Under Harris recurrence of $X$ with invariant probability $\mu$, we determine 'explicite' finite dimensional distributions of the invariant probability $m$ of $\mathbb{X}$ from $\mu$ and the semigroup, and have strong laws of large numbers (\cite{HK-10a}, theorem 2.1) for a large class of functionals 
\beqq\label{myfunctionals}
A = (A_t)_{t\ge 0} \quad,\quad A_t \;=\; \int_0^t F(s,\xi_s)\, \La_T(ds)  
\eeqq
with functions $F(\cdot,\cdot)$ which are $T$-periodic in the first argument ($F(s,x) = F( i_T(s),x )$ for all $s$, $x$) and $\si$-finite measures $\La_T$ on $(\bbr,\calb(\bbr))$ which are $T$-periodic (i.e.\ $\La_T(B)=\La_T(B-kT)$ for all $B$, $k$). The class of functionals (\ref{myfunctionals}) is essentially larger than the class of additive functionals of the process $(\xi_t)_{t\ge 0}$;  
under positive Harris recurrence of the grid chain $X$,  limits   
\beqq\label{mylimit}
\lim_{t\to\infty} \frac1t A_t \;\;=\;\; \frac1T \int_0^T\!\!\!\int_E  \La_T(ds) [\mu P_{0,s}](dy) F(s,y) \quad\mbox{a.s.}
\eeqq
for functionals (\ref{myfunctionals}) will exist provided the integral in the limit is well defined. 

The remaining task is to obtain verifyable criteria which establish:   

{\bf (H)}:\quad the grid chain $X = (\xi_{kT})_{k\in\bbn_0}$ is positive Harris recurrent.  

This is the topic of the present note. Our main result (theorem 1.1) will be stated in section 1; the proof will be given in 2.3 in  section 2. In section 2, we will first review a general route to positive Harris recurrence in terms of Lyapunov functions, following Meyn and Tweedie \cite{MT-92}, then consider SDE with bounded smooth coefficients and uniform ellipticity where classical heat kernel bounds for the transition probabilites, due to Aronson \cite{A-67}, can be used, and then consider SDE with smooth coefficients having linear order of growth, and where the diffusion matrix may have eigenvalues equal to $0$: here we will use non-classical lower bounds on transition probabilities of possibly degenerate diffusions which are due to Bally \cite{B-06}. Some examples are included in section 2. For diffusion processes, we thus dispose of simple conditions in terms of drift and diffusion coefficient which imply positive Harris recurrence {\bf (H)} of the grid chain, thus of the segment chain, and thus of strong laws of large numbers (\ref{mylimit})+(\ref{myfunctionals}) in the process $\xi$.

\section*{1. Main Result}

Consider a $d$-dimensional SDE with $T$-periodic drift 
\beqq\label{process-SDE}
d\xi_t \;=\; b(t,\xi_t)\, dt \;+\; \si(\xi_t)\, dW_t     \quad,\quad   b(t,x) \;=\;  b(i_T(t),x)  
\eeqq
driven by $m$-dimensional Brownian motion with $m\ge d$. We write $\lambda_* (x) \ge 0$ for the smallest and $\lambda^* (x)$ for the largest eigenvalue of $a(x)=\si(x)\si^{\!\top}(x)$,  $\,x\in\bbr^d$.  Let $\sigma_i $ denote the $i-$th line of the $d \times m-$matrix $\sigma $, and $ b = (b^1, \ldots, b^d ).$  We require that the coefficients $b(\cdot,\cdot)$ and $\sigma(\cdot)$ of equation (\ref{process-SDE}) are $d+2$ times differentiable in the space variable with bounded derivatives 
\beam\label{ballycondition-3}
\max_{ 1 \le | \alpha | \le d +2} \max_{i,j} \left( | D^\alpha \sigma_{i,j} (x) | + | D^\alpha b^i (t,x) | \right) \;\le\;  C_0 &&\mbox{for all $0\le t\le T$, all $x\in\bbr^d$}    
\eeam
(which clearly implies Lipschitz and linear growth conditions). Take  $\La$ in (\ref{semigroupTperiodic}) the Lebesgue measure on $\bbr^d$. Then with strong reference to Meyn and Tweedie (\cite{MT-92}, theorem 4.6) and Bally (\cite{B-06}, theorem 24), see section 2 below, we have  \\

{\bf 1.1 Theorem: } {For a $d$-dimensional SDE with $T$-periodic drift (\ref{process-SDE}) which satisfies (\ref{ballycondition-3}), the following conditions are sufficient for  {\bf (H)}:  \\
1. Assume there is some compact set  $\tilde K\subset\bbr^d$ such that 
\beqq\label{lyapunovcondition-5}
2\,x^{\!\top} b(s,x) \;+\; {\rm tr}( a(x) )  \;<\; -\vep   \quad\mbox{on}\; [0,T]{\times}\{ x\in \tilde K^c \}
\eeqq
for some $\vep>0.$ Write 
$$ C = \sup_{ s \le T, x \in \tilde K} |2\,x^{\!\top} b(s,x) \;+\; {\rm tr}( a(x) ) | $$
and let $\tilde R = \left( \sup_{ x \in \tilde K} | x | \right) \vee e  .$ \\
2. Choose $R > \tilde R$ minimal such that 
\beqq\label{minimal} 2 \exp \left(  - \frac{ \log R - \log \tilde R - (  [ (2d)^{3/2} (m+1)^{1/2} C_0 + 4 d^3 (m+1) C_0^2] T )^2 }{ 16 d^3(m+1) C_0^2 T } \right)
\le \frac{\varepsilon}{2 (C + \varepsilon)}  ,
\eeqq
where $C_0$ is the constant of (\ref{ballycondition-3}). Let 
$$ K = \{ x : |x| \le R\} $$
and suppose that
\beqq\label{eq:nondeg}
\lambda_* (x) > 0  \;\; \mbox{for all } x \in K    \;. 
\eeqq
Under these conditions, $ K  $ is a 'small' set and $\nu:=\La(\cdot\cap K)$ a 'small' measure for the kernel $P_{0,T}$ 
$$
P_{0,T}(x,dy) \;\ge \; \al\, 1_{ K}(x)\, \nu(dy) \quad,\quad x,y \in \bbr^d
$$
in the sense of Nummelin (\cite{N-78}, \cite{N-85}), for suitable $\al>0$. }

\vskip0.8cm
Note that there are no ellipticity conditions, the matrix $a(x)=\si(x)\si^{\!\top}(x)$ may have eigenvalues $0$ at points $x$ outside the compact $K$, and the transition density $p_{0,T}(x,y)$ is not necessarily jointly continuous in  $(x,y)$. The proof of theorem 1.1 will be given in 2.3 below; see also remark 2.4.\\

{\bf 1.2 Example: } Consider in theorem 1.1 a drift coefficient containing a deterministic periodic signal 
$$
b(t,x) = S(t) + \hat b(x)
$$
where $t\to S(t)$ is a piecewise continuous $T-$periodic function, and $\hat b(\cdot)$ is Lipschitz. With decomposition $f=f^+-f^-$ in positive and negative part, put 
$$
G_S(x) \;:=\; 2\, \sum_{i=1}^d \left(  x_i^- \max S_i^- +  x_i^+ \max S_i^+ \right) \;.
$$ 
Then the following simple condition 
\beqq\label{easycondition}
2 x^\top \hat b(x) + G_S(x)  + tr(a(x))  \;\;<\;\; -\vep   \quad \mbox{for all $x\notin \tilde K$}  
\eeqq
implies (\ref{lyapunovcondition-5}). Note that $G_S(\cdot)$ in (\ref{easycondition}) remains invariant under scaling of the signal $S(\cdot)$ which is relevant for 'frequency modulation' type problems (\cite{HK-11}, see also Ibragimov and Khasminskii \cite{IH-81} p.~209 in the special case $\hat b(\cdot)\equiv 0$ and $\si(\cdot)\equiv 1$).  \halmos

\section*{2. Proof of theorem 1.1, and some related results}

\subsection*{2.1 Lyapunov conditions}

We restart in the general setting of a strong Markov process $\xi$ taking values in $(E,\cale)$, with $T$-periodic semigroup (\ref{semigroupTperiodic}) as in the introduction, and ask for sufficient conditions implying {\bf (H)} (positive Harris recurrence of the grid chain). The following theorem A, due to Meyn and Tweedie \cite{MT-92}, gives an answer which is by now classical in the literature on Markov processes (it exists in several variants, with different types of Lyapunov conditions in place of (\ref{lyapunovcondition-1}), see e.g.\ \cite{HM-08}). Below, a measurable function $V:E\to[0,\infty)$ is termed 'norm-like' if it is bounded on compacts, and behaves outside large compacts as a strictly increasing function of the distance away from $0$. \\

{\bf Theorem A: } (essentially Meyn and Tweedie \cite{MT-92}, theorem 4.6) Assume that there exists some norm-like function~$V$ and some compact $K\subset E$ such that 
\beam\label{lyapunovcondition-0}
P_{0,T} V \quad\mbox{is bounded on}\;\; K \;, \\
\label{lyapunovcondition-1}
P_{0,T} V \;\le\; V - \vep \quad\mbox{on}\;\; K^c \;, 
\eeam
and such that some lower bounds for the $\Lambda$-density $p_{0,T}(\cdot,\cdot)$ on $K\times K$ establish   
\beqq\label{irreducibilitycondition-2} 
\inf_{x,y\in K}\, p_{0,T}(x,y)\;\;>\;\; 0  \;. 
\eeqq 
Then {\bf (H)} holds. Moreover,  $K$ is a 'small' set and $\nu=\Lambda(\cdot\cap K)$ a 'small' measure for the kernel $P_{0,T}$ 
\beqq\label{nummelincondition-1}
P_{0,T}(x,dy) \;\ge\; \al\, 1_K(x)\, \nu(dy)  \quad,\quad x,y\in E
\eeqq
for suitable $\al>0$. \\

{\bf Idea of proof: } 1) $(V(\xi_{kT}))_k$ evolves as a nonnegative supermartingale as long as $\xi_{kT}$ is in $K^c$, as a consequence of condition (\ref{lyapunovcondition-1}), hence a.s.\ the grid chain $(\xi_{kT})_k$ with starting point $x\in K^c$ will enter the compact $K$ in finite time.  \\
2) By condition (\ref{irreducibilitycondition-2}), we have $\inf\limits_{x,y\in K}p_{0,T}(x,y)\;$ strictly positive which allows for minorization (\ref{nummelincondition-1}).\\
3) Nummelin splitting, from (\ref{nummelincondition-1}) combined with step 1), shows Harris recurrence of the grid chain  $(\xi_{kT})_k$, with a unique (up to constant multiples) invariant measure $\mu$ which is equivalent to the maximal irreducibility measure $\nu \sum_{k=0}^\infty 2^{-k}P_{0,kT}$ (see Nummelin \cite{N-78}, \cite{N-85}).\\
4) Now we refer to the Meyn-Tweedie result only for {\em positive} Harris recurrence, i.e. for the argument \cite{MT-92} p.\ $555_{6-4}$ and p.\ $557_{13}$. At this last step, we need (\ref{lyapunovcondition-0}) which on $K$ yields a bound for the expected duration  of excursions away from $K$; see also \cite{N-85} p.\ $78^{11-13}$. \halmos\\

Obviously, condition (\ref{irreducibilitycondition-2}) in theorem A is implied by the  stronger condition 
\beqq\label{irreducibilitycondition-1} 
(x,y) \;\lra\; p_{0,T}(x,y)  \quad\mbox{is strictly positive and continuous} \; .
\eeqq 
We illustrate theorem A by two examples where (\ref{irreducibilitycondition-1}) is at hand, where the drift has the structure of example 1.2, with different choices of a norm-like function $V(\cdot)$. \\

{\bf 2.1 Examples: } We consider $E=\bbr$, $\Lambda=\lala$ Lebesgue measure, $\xi$ an Ornstein Uhlenbeck type process.   
The driving semimartingale is either one-dimensional Brownian motion $W$ or a one-dimensional L\'evy process $Z$ without Gaussian component. 
 
a)~ Choose $V(y)=y^2$.  For some constants  $\gamma>0$, $\si>0$, let 
$$
d\xi_t \;=\; \left(\,S(t) -\gamma\, \xi_t\,\right)\, dt \;+\; \si\, dW_t  \quad,\quad t\ge 0 \;. 
$$
Here all transition probabilities $P_{s,t}(\cdot,\cdot)$ are normal laws (\cite{HK-10a}, example 2.3), and in particular 
$$
P_{0,T}(x,\cdot) \;=\; \caln\left( x e^{-\gamma T} + \int_0^T e^{-\gamma v}  S(T-v)  dv \,,\; \frac{1-e^{-2 \gamma T}}{2\gamma}\;\si^2 \right) \;. 
$$
Hence (\ref{irreducibilitycondition-1}) is obvious, thus (\ref{irreducibilitycondition-2}). The elementary relation between variance and second moments of a random variable gives 
$$
P_{0,T} V (x) \;=\; \left[\frac{1-e^{-2 \gamma T}}{2\gamma}\;\si^2\right] + \left[x e^{-\gamma T} + \int_0^T e^{-\gamma v}  S(T-v)  dv\right]^2 
$$
which is bounded on compacts, and such that for $|x|$ tending to $\infty$ and suitably large $M<\infty$
$$
P_{0,T} V (x) \;\;=\;\; x^2\, e^{-2 \gamma T} \;+\; O(|x|)  \quad\le\quad V(x) \;-\; \vep  \quad\mbox{for $|x|>M$}  \;. 
$$  
This gives (\ref{lyapunovcondition-0})+(\ref{lyapunovcondition-1}). By theorem A, we have {\bf (H)} and (\ref{nummelincondition-1}). 
With $T$-periodic means $s\to M(s)$
\beqq\label{functionM}
M(s) \;:=\;  \int_0^\infty e^{-\gamma v} S(s-v) dv \;=\; \int_0^T \frac{e^{-\gamma v}}{1-e^{-\gamma T}}\, S(s-v)\, dv  
\eeqq
the invariant probability $\mu$ of the grid chain $X$  and the measures $\mu P_{0,s}$ in (\ref{mylimit}) are given by  
$$
\mu \;=\; \caln\left( M(0)\,,\, \frac{\si^2}{2\gamma}\, \right) \quad,\quad 
\mu P_{0,s} \;=\; \caln\left(\, M(s) \,,\,  \frac{\si^2}{2\gamma}\,\right) \;,\; 0\le s\le T \;,
$$
see \cite{HK-10a}. In this example, the $T$-segment chain $\mathbb{X}$ is positive Harris taking values in $(\mathbb{C}_T,\calc_T)$. 

b)~ We take $V(x) = |x|$ and consider for constant $\gamma>0$
\beqq\label{OUjumpdiffusion}
d\xi_t \;=\; ( S(t)-\gamma \xi_t )\, dt \;+\; dZ_t  
\eeqq
where $Z$ is a one-dimensional L\'evy process with L\'evy triplet $( 0, 0, \nu)$ whose L\'evy measure satisfies 
\beqq\label{integrierbarkeitsbedingung}
\int_{\{ |x|>1 \}} |x| \nu (dx) \;<\; \infty  
\eeqq
and 
\begin{equation}\label{eq:bk}
[\varepsilon^2 \ln \frac1\varepsilon ]^{-1} \int ( x^2 \wedge \varepsilon^2 ) \nu (dx) \;\to\;  + \infty \quad\mbox{ as  } \varepsilon \to 0 + \;.
\end{equation} 
The condition (\ref{eq:bk}) was given by Bodnarchuk and Kulik (\cite{BK-08}, condition (iii) in theorem 1). In the special case $S\equiv 0$, taking as starting point $x=0$ at time $s$, the process $\xi$ possesses a transition density $p^0_{s,t} (0, y) $ which is $C^\infty_b$ in $y\;$ (\cite{BK-08}, theorem 1). By explicit representation 
\beao
\xi_t &=& e^{-\gamma(t-s)}\, \xi_s \;+\; \int_s^t e^{-\gamma(t-v)}\, S(v)\, dv \\
&+&  \int_s^t \int_{\{ |x|>1 \}} e^{-\gamma(t-v)}\, x\, [dv\, \nu(dx)]  \;\;+\;\; \int_s^t \int_{\{ |x|\le 1 \}} e^{-\gamma(t-v)}\, x\, (\mu(dv,dx)-dv\, \nu(dx)) 
\eeao
of the solution to (\ref{OUjumpdiffusion}) (as in \cite{BK-08}) for general signal $S$ and  pairs $s<t,$ the transition density of the process (\ref{OUjumpdiffusion}) can be written as 
$$
p_{s,t}(x,y) \;=\; p^0_{s,t}\left(0, \, \; y - e^{-\gamma(t-s)}\, x - \int_s^t e^{-\gamma(t-v)}\, S(v)\, dv\, \right) \;.
$$
This implies joint continuity of $(x,y) \to p_{s,t}(x,y)$ which is (\ref{irreducibilitycondition-1}), and thus (\ref{irreducibilitycondition-2}) in theorem A for arbitrary compacts $K$. Condition (\ref{integrierbarkeitsbedingung}) guarantees $P_{0,t}V(x) = E_x ( | \xi_t| ) < \infty $ for all $t$ and $x$, and (again from the explicit representation of the solution) implies condition (\ref{lyapunovcondition-0}) of theorem A. Also 
\begin{eqnarray*}
 P_{0,T} V (x) &=& \int p^0_{0,T} \left(0, y - e^{ - \gamma T } x - \int_0^T e^{-\gamma(T-s)}\, S(s)\, ds\right) |y| dy\\ 
&= &\int p_{0,T}^0 (0, y) | y + e^{ - \gamma T } x + \int_0^T e^{-\gamma(T-s)}\, S(s)\, ds | dy \\
&\le &\int p_{0,T}^0 (0, y)\left( e^{-\gamma T} |x| + |y| + \int_0^T e^{-\gamma(T-s)}\, |S(s)| \, ds\right) dy \\
&= & V(x) + \int p_{0,T}^0 (0, y) | y | dy  + \int_0^T e^{-\gamma(T-s)}\, |S(s)| \, ds + [ e^{ - \gamma T } - 1 ]| x|   \\
&\le& V(x) - \varepsilon 
\end{eqnarray*}
provided $x \in K^c,$ for $K$ an appropriate compact set. This gives condition (\ref{lyapunovcondition-1}) of theorem A. Hence by theorem A, we have {\bf (H)} and (\ref{nummelincondition-1}). Here, the $T$-segment chain $\mathbb{X}$ is positive Harris taking values in $(\mathbb{D}_T,\cald_T)$. 
We specify the invariant probability $\mu$ for the grid chain $X$ together with the laws $\mu P_{0,s}$ arising in (\ref{mylimit}). Extend Poisson random measure $\mu(ds,dx)$ with intensity $ds\, \nu(dx)$ from $(0,\infty)\times(\bbr{\setminus}\{0\})$ to $\bbr\times(\bbr{\setminus}\{0\})$ and consider 
$$
U \;:=\; \int_{-\infty}^0 \int_{\{|x|>1\}} [x e^{ - \gamma |v| }]\, dv\, \nu(dx) 
\;+\; \int_{-\infty}^0 \int_{\{|x|\le 1\}} [x e^{ - \gamma |v| }]\, ( \mu(dv,dx) - dv\, \nu(dx) ) \;. 
$$ 
From the explicit representation of the solution to (\ref{OUjumpdiffusion}), we see that the law 
$$
\mu \;=\; \call\left(\, M(0) + U \,\right) 
$$
is invariant for $P_{0,T}$, where $M(\cdot)$ is the $T$-periodic function defined in (\ref{functionM}). By Harris recurrence {\bf (H)}, $\mu$ is the unique invariant probability for the grid chain $X$.  Similarly,    
$$
\mu P_{0,s} \;=\; \call\left(\, M(s) + U \,\right) \;,\; 0\le s\le T \;,
$$
thus the limit in (\ref{mylimit}) is explicit. Our integrability conditions on BV part and martingale part of  
$$
Z_t \;=\; \int_0^t \int_{\{ |x|>1 \}} x\, [dv\, \nu(dx)]  \;\;+\;\; \int_0^t \int_{\{ |x|\le 1 \}} x\, (\mu(dv,dx)-dv\, \nu(dx)) 
$$
give 
$$
\int_0^\infty e^{ - \gamma v }\, dZ_v  \;=\;  \int_0^\infty Z_v\, \gamma e^{ - \gamma v } dv 
$$
(see \cite{DM-80}, (19.2) in ch.\ VIII), and thus the interpretation 
$$
\call \left(\, U \,\right)\;=\;  \call \left(\, Z_\tau \, \right)  
$$ 
where $Z_\tau$ is the process $Z$ evaluated at an independent exponential time  $\tau$ with parameter $\gamma$. 
 
c)~Note that for a symmetric stable $Z$ having $ \nu (dx) = c |x|^{ - \alpha - 1 },$ $0 < \alpha < 2,$ both conditions (\ref{integrierbarkeitsbedingung}) and (\ref{eq:bk}) are satisfied simultaneously iff $ 1 < \alpha < 2 $. \halmos\\

\subsection*{2.2 A classical ansatz for SDE: heat kernel bounds}

Now we restrict the setting to $d$-dimensional SDE with $T$-periodic drift (\ref{process-SDE})
$$
d\xi_t \;=\; b(t,\xi_t)\, dt \;+\; \si(\xi_t)\, dW_t     \quad,\quad   b(t,x) \;=\;  b(i_T(t),x)  
$$
driven by $m$-dimensional Brownian motion ($m\ge d$, $\,a(x)=\si(x)\si^{\!\top}(x)$). According to Aronson \cite{A-67} we suppose uniform ellipticity 
\beqq\label{aronsoncondition-1}
\zeta^{\!\top}\, a(x)\, \zeta \;\;\ge\;\; cst \cdot |\zeta|^2   \quad\mbox{for all $\zeta,x$ in $\bbr^d$} 
\eeqq
and smoothness of the coefficients of equation (\ref{process-SDE}) as follows: for some $\delta>0$ small, 
\beam 
&&(t,x)\to b(t,x) \;\mbox{and}\; x\to a(x)\;\mbox{are smooth and bounded on}\; (0,T+\delta){\times}\bbr^d
\label{aronsoncondition-2} \\
&&\mbox{all partial derivatives of order $1$ of the functions in (\ref{aronsoncondition-2}) are bounded on}\; (0,T+\delta){\times}\bbr^d \;. 
\label{aronsoncondition-3}
\eeam
Note that condition (\ref{aronsoncondition-3}) arises from Aronson \cite{A-67} when his equation (1) is written as a particular case of his equation (3). By \cite{A-67}, theorem 1, we have 'heat kernel bounds' on $(0,T+\delta)\times\bbr^d$: \\
for all $x,y$ in $\bbr^d$ and all $0<s<T{+}\delta$, 
\beqq\label{heat-kernel-bounds-1}
\kappa_1  \left(\frac{1}{2\pi\, s\,\gamma_1^2}\right)^{d/2} \!\exp\{ -\frac12 \frac{|y-x|^2}{s\,\gamma_1^2} \} 
\;\;\le\;\; p_{0,s}(x,y) \;\;\le\;\; 
\kappa_2  \left(\frac{1}{2\pi\, s\,\gamma_2^2}\right)^{d/2} \!\exp\{ -\frac12 \frac{|y-x|^2}{s\,\gamma_2^2} \} 
\eeqq
for constants $\kappa_1,\kappa_2,\gamma_1,\gamma_2$ in $(0,\infty)$ which depend on the dimension $d$, the constant in (\ref{aronsoncondition-1}), the time horizon $T{+}\delta$, and the bounds in (\ref{aronsoncondition-2})+(\ref{aronsoncondition-3}) on the coefficients of equation (\ref{process-SDE}). 
Notice that similar bounds have been obtained in \cite{KS-85} and \cite{NS-91}. \\

{\bf 2.2 Theorem: } Under Aronson's conditions (\ref{aronsoncondition-1})--(\ref{aronsoncondition-3}) assume in addition that for suitably large $\tilde R$
\beqq\label{lyapunovcondition-3}
2\,x^{\!\top} b(s,x) \;+\; {\rm tr}( a(x) )  \;<\; -\vep   \quad\mbox{on}\; [0,T]{\times}\{ |x|>\tilde R \}  \;. 
\eeqq
Then (\ref{lyapunovcondition-0})+(\ref{lyapunovcondition-1})+(\ref{irreducibilitycondition-2}) in theorem~A hold, hence {\bf (H)} and (\ref{nummelincondition-1}). \\

{\bf Proof: } Consider $V(y)=|y|^2$. Write $\call_s V (x) = 2\,x^{\!\top} b(s,x) + \sum_{i=1}^d  a^{i,i}(x)$ which by our assumptions on the coefficients in equation (\ref{process-SDE}) is bounded on $[0,\infty){\times}\{ |x|\le \tilde R \}$. Similarly, by $T$-periodicity of $b(\cdot,x)$, (\ref{lyapunovcondition-3}) holds also on $[0,\infty){\times}\{ |x|> \tilde R \}$. Then Ito formula shows 
\beao
V(\xi_t) - V(\xi_0)  &=&   \int_0^t 2\, \xi_s^\top \si(\xi_s)\, dW_s
\;+\; \int_0^t 1_{\{|\xi_s|\le \tilde R\}} (\call_s V)(\xi_s)\, ds  \;+\; \int_0^t 1_{\{|\xi_s|> \tilde R\}} (\call_s V)(\xi_s)\, ds\\
&\le&  \int_0^t 2\, \xi_s^\top \si(\xi_s)\, dW_s  \;+\; C \int_0^t 1_{\{|\xi_s|\le \tilde R\}}\, ds \;-\; \vep\, \int_0^t 1_{\{|\xi_s|> \tilde R\}}\, ds  \\
&=&   \int_0^t 2\, \xi_s^\top \si(\xi_s)\, dW_s  \;-\; \vep\, t  \;+\; (C+\vep) \int_0^t 1_{\{\,|\xi_s|\le \tilde R\,\}}\, ds \;,
\eeao
for suitable $C<\infty$. With stopping times $S_N:=\inf\{t>0:|\xi_t|>N\}$ and stopped processes $\xi^{S_n}=(\xi_{t\wedge S_N})_{t\ge 0}$, the local martingale in the Ito formula is a martingale up to time $S_N$, thus by Fatou  
\beqq\label{prae-lyapunov-condition}
\left( P_{0,T}V - V \right)(x)  \;\;\le\;\;  -\; \vep\, T \;\;+\;\; (C+\vep) \int_0^T P_x(\,|\xi_s|\le \tilde R \,)\, ds 
\eeqq
for $t=T$. Aronson's upper bound in (\ref{heat-kernel-bounds-1}) shows 
$$
\mbox{for every $0<s\le T$ :}\quad
P_x(\,|\xi_s|\le \tilde R \,) \;\lra\; 0 \quad\mbox{as}\; |x|\to\infty \;. 
$$
Thus dominated convergence in (\ref{prae-lyapunov-condition}) shows that condition (\ref{lyapunovcondition-1}) in theorem A holds, with $K = \{ |x| \le  R \},$ for some suitable $R  > \tilde R .$ $\,V(\cdot)$ being bounded on compacts, (\ref{prae-lyapunov-condition}) shows that the function $P_{0,T}V$ is bounded on $K$: this is condition (\ref{lyapunovcondition-0}) in theorem A. Condition (\ref{irreducibilitycondition-2}) in theorem A is an immediate consequence of Aronson's lower bound with $s=T$ in (\ref{heat-kernel-bounds-1}).  Hence theorem~A applies and gives the result. \halmos

\subsection*{2.3 More general SDE: non-classical lower bounds for the transition density}

In this subsection we prove theorem 1.1 under condition (\ref{ballycondition-3}) through lower bounds on the transition probability of a $d$-dimensional diffusion (\ref{process-SDE}) with $T$-periodic drift. The following theorem, due to Bally \cite{B-06}, establishes a lower bound on the transition probabilities $p_{0,T}(x,y)$ on suitable compact and convex sets in $\bbr^d$. We do not need (uniform) ellipticity, and do not have to assume boundedness of the coefficients of equation (\ref{process-SDE}). All notations are as in section 1.\\

{\bf Theorem B: } (Corollary to Bally \cite{B-06}, theorem 24) Consider equation (\ref{process-SDE}) with coefficients such that (\ref{ballycondition-3}) holds, and some compact and convex set $K\subset \bbr^d$ which satisfies (\ref{eq:nondeg})
$$
\lambda_* (x) > 0  \;\; \mbox{for all $x \in K$} \;. 
$$
Then there exists a strictly positive constant $C = C( K, d,m, C_0, \lambda_*, \| \lambda^* \|_{\infty, K} ) $ such that 
\begin{equation}\label{eq:lowerbound}
p_{0,T} ( x, y) \;\geq\;  C  \quad\mbox{for all $x,y \in K$} \;. 
\end{equation}

\vskip0.8cm
{\bf Proof:} We  check the set of conditions in \cite{B-06}, theorem 24. In view of a compact set to be used in (\ref{eq:nondeg}), it is sufficient to consider $N(x)$ defined by $N^2(x) = 1 + |x|^2$. Then assumption 
(\ref{ballycondition-3}) implies Lipschitz  
\beqq\label{ballycondition-2}
\quad\max_{ 1 \le i \le d} \left( \| \sigma_i (x) - \sigma_i (y) \| + | b^i (t,x) - b^i (t,y) | \right)  \;\;\le\;\;  C_0 \sqrt{ (m+1) d} \, \, | x - y |  \quad\mbox{for all $0\le t\le T$,  $x,y\in\bbr^d$ }  
\eeqq
and linear growth  conditions
\beqq\label{ballycondition-1}
\max_{1 \le i \le d} \left( \|\sigma_i(x)\| + |b^i(t,x)| \right) \;\;\le\;\;  C_0 \sqrt{ (m+1) d}  \, N(x)  \quad  \mbox{for all $0\le t\le T$,  $x\in\bbr^d$} \;,
\eeqq
thus conditions A of \cite{B-06}. Fix $x_0 , y \in K .$ Consider the line 
$ 
x_t \;=\; x_0 + \frac{t}{T} (y-x_0) ,$ $ t \in [0, T] \, 
$
from $x_0$ to $y.$ Note that by convexity of $K,$ $x_t \in K $ for all $ t\in [0, T],$ and hence
$a^{-1} (x_t)$ exists for all $t \in [0, T]$, by (\ref{eq:nondeg}) above. Hence we can define the control 
$$ 
\phi_t \;=\; \sigma^{\top} (x_t)\, a^{-1} (x_t)\, \partial_t x_t  \;\;,\;\; \partial_t x_t = \frac1T (y -x_0) 
$$ 
and have with this definition a differentiable path $x_t = x_t^\phi$ where $x_t^\phi $ is solution to the equation 
$$ 
dx_t^\phi \;=\; \sigma( x_t^\phi)\, \phi_t\, dt \;.
$$
We  put 
$ 
\varrho (x) \;=\; \frac{\sqrt{ \lambda_* (x)}}{   N(x)} \;.
$
Since $K$ is compact and $N(x)$ bounded on compacts, condition (\ref{eq:nondeg}) implies that for suitable constants $ \mu \geq 1 $ and $\chi > 0, $ 
$$ \varrho (x) \geq \frac{1}{\mu} \mbox{  and  } \lambda_* (x) \geq \frac{1}{\chi^2} \mbox{  for all } x \in K . 
$$ 
Since $x_t = x_t^\phi  \in K$ for all $t \in [0, T],$  this implies  
$$
\rho(x_t^\phi) \;\ge\; \frac1\mu \;\;,\;\;  \sqrt{ \la_*(x_t^\phi) } \;\ge\; \frac1\chi 
$$ 
which is the condition (\cite{B-06}, p.\ $2435_{7}$). Next, since  
$
| \phi_t |^2 \;=\; | \phi_t^{\top}\phi_t |  \;=\; \| a^{-1} (x_t) \| \cdot | \frac1T ( y - x_0 ) |^2 ,
$
we have for any pair $s,t$ in $[0,T]$ 
$$ 
\frac{| \phi_t |}{| \phi_s| } \;=\; \frac{\| a^{-1} (x_t) \|^{1/2} }{\| a^{-1} (x_s) \|^{1/2} } \;\le\; \eta \;<\; \infty 
\quad,\quad  \eta^2 := \sup_{ x, y \in K}  \frac{\| a^{-1} (x) \| }{\| a^{-1} (y) \| } 
$$ 
since $x_t, x_s \in K$ and $K$ is compact, by (\ref{eq:nondeg}) and boundedness of $\la^*(x)$ on $K$. The same argument yields that there exists a constant $\nu$ such that
$ 
| \phi_t | \;\le\; \nu \;, \;\; \forall 0\le t \le T \;.
$
Hence we have checked (\cite{B-06}, p.\ $2435_{6}$). The conditions preceding (\cite{B-06}, theorem 24) are now satisfied with the set of parameters $\theta = ( \mu , \chi, \nu , \eta , h ),$ for any fixed $ h \in ] 0, T[ .$  
It remains to verify the condition (\cite{B-06}, p.\ $2435_{1}$) in the beginning of this theorem. To check this, we have to upper-bound the quantity $d_\theta (x_0, y )$ defined in (\cite{B-06}, p.\ $2435_{4-3}$) where by definition (cf.\  \cite{B-06}, p.\ $2409_8$)
$$
d^2_\theta (x_0, y ) \;\le\;  \int_0^T | \phi_t |^2\, dt 
$$
for the above control $\phi_t. $ But, as above 
$$ 
| \phi_t |^2 \;\le\; \frac{1}{T^2}\, \| a^{ -1} (x_t) \|\, | y-x_0 |^2 \;\le\; \frac{C}{T^2} 
$$
for a suitable constant $C$ where the bound does not depend on $x_0,y$ in the compact $K$, and thus 
$$ 
d_\theta (x_0, y) \;\le\; C\, \frac{1}{\sqrt{T}} 
$$ 
simultaneously for all $x_0, y \in K .$ Thus all conditions of  (\cite{B-06}, theorem 24) are checked. Now (\ref{eq:lowerbound}) follows from the strictly positive lower bound  (\cite{B-06}, formula (21)) on $p_{0,T}(x_0,y)$ simultaneously for all $x_0, y \in K .$ This concludes the proof. \halmos\\

In the above proof no smoothness of $b(t,x)$ with respect to the time variable is required. Indeed, \cite{B-06} uses conditional Malliavin calculus where derivatives are only taken with respect to the space variables. Using theorem B we can prove theorem 1.1. \\

{\bf 2.3 Proof of theorem 1.1: } 1) {We consider equation} (\ref{process-SDE}) with coefficients satisfying (\ref{ballycondition-3}). We assume that there is some compact set  $\tilde K\subset \bbr^d$ satisfying (\ref{lyapunovcondition-5}) 
$$
2\,x^{\!\top} b(s,x) \;+\; {\rm tr}( a(x) )  \;<\; -\vep   \quad\mbox{on}\; [0,T]{\times}\{ x\in \tilde K^c \}  \;
$$
and some constant $R$  (see (\ref{minimal}))  defined in terms of $\tilde K$, $C_0,$ $m$ and the dimension $d$ such that for $K = \{ x : |x| \le R \}, $ assumption (\ref{eq:nondeg})  
$$
\lambda_* (x) > 0  \;\; \mbox{for all $x \in K$}  
$$
holds. 
 
We put $V(x)=|x|^2$ and proceed on the lines of the proof of theorem 2.2 to obtain equation (\ref{prae-lyapunov-condition}):  
$$
\left( P_{0,T}V - V \right)(x)  \;\;\le\;\;  -\; \vep\, T \;\;+\;\; (C+\vep) \int_0^T P_x(\,\xi_s\in \tilde K \,)\, ds  
$$
for arbitrary $x\in\bbr^d$. {Here, $C$ is a bound on $ {\cal L}_s V (\cdot) $ on $\tilde K.$} 

2) We prove that under the present conditions we have  
\beqq\label{letzterknackpunkt}
\mbox{for every $0<s\le T$ :}\quad
P_x(\, \xi_s\in \tilde K \,) \le  \frac{ \vep }{ 2(C+\vep) } \quad\mbox{for all } x \in K^c  \;.  
\eeqq
Recall the definition of $ \tilde R$ in theorem 1.1: We have $\tilde R\in(e,\infty)$ and $\tilde K\subset B_{\tilde R}(0).$ We now consider a $\calc^2$--function $F:\bbr^d\to(0,\infty)$ with the properties $F(x)=\log(|x|)$ on $\{|x|>\tilde R\}$ and $F(x)\le \log(\tilde R), $ {$|D_i F(x) | \le \frac{2}{\tilde R}, $ $ |D_{i,j} F (x) | \le \frac{4}{\tilde R^2} $} on $B_{\tilde R}(0)$. Thanks to linear growth  (\ref{ballycondition-1}) of the coefficients of equation (\ref{process-SDE}), and to the order of growth of partial derivatives $D_iF$ and $D_{i,j}F$ on  $\{|x|>\tilde R\}$, for $1\le i,j\le d$,  all integrands in the Ito formula 
\beao
F(\xi_t) - F(\xi_0) &=& \sum_{i=1}^d \int_0^t D_iF(\xi_s)\, b^i(s,\xi_s)\, ds  
\;+\; \frac12 \sum_{i,j=1}^d \int_0^t D_{i,j}F(\xi_s)\, a_{i,j}(\xi_s)\, ds  \;+\; M_t\\
M_t   &:=& \sum_{i=1}^d \int_0^t D_iF(\xi_s)\, \si_i(\xi_s)\, dW_s 
\eeao
are bounded {by a constant depending only on the dimension $d$ and on the linear growth bound $C_0 \sqrt{ (m+1) d } $ of (\ref{ballycondition-1})}. 
More precisely,  
$$
|  \sum_{i=1}^d D_iF(x )\, b^i(s,x )| \le   (2d)^{3/2} (m+1)^{1/2} C_0,$$
where we have used that $ N(x) / |x| \le \sqrt{2} $ for $ |x| \geq 1 $ and $|D_i F(x) | \le \frac{2}{\tilde R} $ on $B_{\tilde R}(0).$
The same argument gives
$$ |  \sum_{i=1}^d D_iF(x )\, \si_i(x )| \le (2d)^{3/2} (m+1)^{1/2} C_0 . $$
Similarly, 
$$
| \frac12 \sum_{i,j=1}^d D_{i,j}F(x )\, a_{i,j} (x )| \le 4 d^3 (m+1)  C_0^2 .
$$
By choice of $F(\cdot)$ and $\tilde R$, we can write for $0\le s\le T$ and for $|x| > R$
\beao
P_x(\, \xi_s\in \tilde K \,) &\le& P_x(\,F(\xi_s)\le F(\tilde R)\,) \;\;=\;\; P_x\left(\, F(x)-F(\xi_s) \;\ge\; F(x)-F(\tilde R) \,\right)\\
&&\le\quad   P_x\left(\, \sup\limits_{0\le s\le T}|M_s| \;\ge\; [F(x)-F(\tilde R)-C_2T] \,\right)
\eeao
where $C_2 = (2d)^{3/2} (m+1)^{1/2} C_0 + 4 d^3 (m+1) C_0^2$ is the bound on the integrands of the BV terms in the Ito formula above. $M$ has angle bracket 
$$\langle M \rangle_t = \int_0^t  \sum_{i,j} D_i F (\xi_s) D_j F( \xi_s) a_{i,j} ( \xi_s) ds \le  8 d^3 (m+1) C_0^2 t .$$
So we apply a Bernstein type inequality (see for example inequality (1.5) in Dzhaparidze and van Zanten \cite{DVZ}) and get for all $|x|>R$
\begin{eqnarray*}
&&P_x\left(\, \sup\limits_{0\le s\le T}|M_s| \;\ge\; [F(x)-F(\tilde R)-C_2T] \,\right)\\
&&\le 2 \exp \left(  - \frac{(F(x)-F(\tilde R)-C_2T)^2}{2  (8 d^3 (m+1) C_0^2 T)  }\right) \\
&&= 2 \exp \left(  - \frac{ ( \log |x| - \log \tilde R - [ (2d)^{3/2} (m+1)^{1/2} C_0 + 4 d^3 (m+1) C_0^2 ] T )^2 }{ 16 d^3 (m+1) C_0^2 T } \right) .
\end{eqnarray*}
Finally, by the choice of $R$ in (\ref{minimal}), the last expression is upper bounded by 
$$
P_x\left(\, \sup\limits_{0\le s\le T}|M_s| \;\ge\; [F(x)-F(\tilde R)-C_2T] \,\right) \le  \frac{ \vep }{ 2(C+\vep) }  , \quad \mbox{ for all } |x| > R .
$$
Putting the last inequalities together, we have proved (\ref{letzterknackpunkt}). Combining (\ref{letzterknackpunkt}) with equation (\ref{prae-lyapunov-condition}) in step 1), we have established condition (\ref{lyapunovcondition-1}) in theorem A, with $K = \{ x : |x| \le  R \}.$ 

3) By equation (\ref{prae-lyapunov-condition}) in step 1) above, condition (\ref{lyapunovcondition-0}) is satisfied on $K.$ Hence, both conditions (\ref{lyapunovcondition-0}) and (\ref{lyapunovcondition-1}) in theorem A are checked. 
Since $K$ is compact and convex and since (\ref{eq:nondeg}) holds by assumption on $K$, we can apply theorem B which establishes (\ref{irreducibilitycondition-2}) in theorem A. Then the assertion of theorem 1.1 follows from theorem A. \halmos\\

{\bf 2.4 Remark: }  {\it Note that the assertion to have a compact and convex set $K$ which is 'small' (as in theorem 1.1) is much weaker than assuming ellipticity for SDE's (assumption (\ref{aronsoncondition-1}) in theorem 2.2). This is illustrated by the following example which is taken from Bally and Kohatsu-Higa \cite{B-KH-10}. 

1) Let $ \xi_t = (\xi_t^1, \xi_t^2) $ be the solution of the following degenerate stochastic differential equation
\beqq\label{BKHprocess}
d \xi_t^1 \;=\; b^1 (\xi_t)\, dt \;+\;  \sigma (\xi_t )\, d W_t \quad,\quad d \xi_t^2 \;=\; b^2 (\xi_t )\, dt \;.  
\eeqq
Here, the driving Brownian motion is supposed to be one-dimensional. 
In this example, ellipticity and even the strong H\"ormander condition fail at every point $x \in \bbr.$ Suppose that the coefficients $\sigma, b^1, b^2$ are five times differentiable, not necessarily bounded, but with bounded derivatives of orders $1,\ldots,5$, and that the following non-degeneracy condition holds true: \\
there exists a constant $c > 0 $ such that for all $x \in \bbr^2,$ 
$$ 
| \sigma (x) | \;\geq\; c \;\;, \;\; | \frac{\partial b^2}{\partial x^1}(x) | \;\geq\; c\; .
$$
Then for all compact sets $K,$ there is a strictly positive constant $c_K $ such that 
\beqq\label{boundforBKHprocess} 
\inf_{x, y \in K}\, p_{0, T} (x,y) \;\geq\; c_K   
\eeqq
in virtue of (\cite{B-KH-10}, section 1 and theorem 17). 

2) This result (\ref{boundforBKHprocess}) remains valid in the time inhomogenous case with $T$-periodic drift 
$$
b^i(s,x) \;=\; b^i(\, i_T(s) \,,\, x \,) \quad,\quad i=1,2
$$
provided that bounds as above for $\frac{\partial b^2}{\partial x^1}(s,x)$ and for derivatives of $b^i(s,x)$ with respect to the space variable of orders $1,\ldots,5$ hold uniformly in $s\in [0,T]$. Here, for the same reason as indicated after the proof of theorem B, no smoothness of $b(s,x)$ with respect to the time variable is required.  

3) Given (\ref{boundforBKHprocess}) for a bivariate process with $T$-periodic drift 
$$
d \xi_t^1 \;=\; b^1 (t,\xi_t)\, dt \;+\;  \sigma (\xi_t )\, d W_t \quad,\quad d \xi_t^2 \;=\; b^2 (t,\xi_t )\, dt \;,  
$$
a sufficient condition for positive Harris recurrence of the grid chain  {\bf (H)} together with (\ref{nummelincondition-1}) is that 
$$
2\,\sum_{i=1,2} x_i\,  b^i(t,x) \;+\;  \si^2(x)   \;<\; -\vep   \quad\mbox{on $[0,T]\times\hat K ^c$} 
$$ 
holds for some compact $\hat K$ in $\bbr^2$. This follows --as a consequence of theorem A-- by exactly the same reasoning which in the proof of theorem 1.1 leads to (\ref{prae-lyapunov-condition})+(\ref{letzterknackpunkt}).}\\

{\bf 2.5 Remark: } {\it Consider in dimension $d=1$ geometric Brownian motion $X = (X_t)_{t\ge 0}$  
$$
X_t \;=\; X_0\, \exp\left\{\, \si\, B_t + (\mu-\frac12 \si^2 )\,t \,\right\}  \;,\; t\ge 0 \quad,\quad X_0\ne 0   
$$  
in case $\mu<-\frac12\si^2$. In theorem B above,  we have (\ref{ballycondition-3}) and $\la_*(x)=\si^2x^2.$ Hence condition (\ref{eq:nondeg}) can be satisfied only on compact intervals $K$ which do not contain the point~$0$. For $\mu<-\frac12\si^2$, condition (\ref{lyapunovcondition-5}) requires compacts which contain a neighbourhood of~$0$. Both requirements are clearly incompatible. For $\mu\le 0$, geometric Brownian motion has paths such that $\lim\limits_{t\to\infty}X_t=0$ almost surely, thus the grid chain $X=( \xi_{kT} )_k$ can not be Harris.}\\

{\bf 2.6 Remark: } {\it In a time homogeneous setting, Veretennikov \cite{V-97} (cf.\ also \cite{V-84}) imposes the condition (\ref{aronsoncondition-2}). He does not require (\ref{eq:nondeg}) but only $\la^*(x)>0$ for all $x\in\bbr^d$. 
Instead of (\ref{lyapunovcondition-5}) he considers 
$$
\la_- := \inf_{x\neq 0} \frac{ x^\top a(x)\; x}{|x|^2} \;\;,\;\; 
\la_+ := \sup_{x\neq 0} \frac{ x^\top a(x)\; x}{|x|^2} \;\;,\;\; 
\wt{\La} := \sup_x\; {\rm tr}(a(x))
$$
and assumes that 
$$
x^\top b(x) \;\le\; -r \;,\; |x|>M \quad,\quad 
\frac32 \la_+ \;<\; r - \frac12 (\wt{\La} - \la_-)  
$$
(\cite{V-97}, inequality (6), and (23) with $r_0>\frac32$); his condition implies in particular that 
$$
 2 x^\top b(x) + {\rm tr}(a(x)) \;\;<\;\; -2 \la_+ \quad\mbox{on $[0,T]\times \{|x|>M\}$} \;. 
$$
}

\vskip0.8cm 
{\bf 2.7 Example:} Consider in dimension $d=1$ a Pearson type diffusion  
$$ 
d X_t \;=\; \theta\, ( S(t) - X_t )\, dt \;+\; \sigma (X_t) \sqrt{ c_0 +  ( X_t - c_1)^2 }\, d W_t 
$$
with constants $\theta > 0$ and $c_0 > 0$, where $t\to S(t) $ as in 1.2 is a deterministic, strictly positive and continuous $T-$periodic signal. We assume $\sigma(\cdot)$ strictly positive and bounded; it may decrease to $0$ as $|x|\to\infty$, and we assume  $\max\left( |\si '|, |\si ''|, |\si '''| \right) (x) \;=\; O(|x|^{-1}) $ for $ |x|\to\infty$. Then (\ref{ballycondition-3}) is satisfied. (\ref{eq:nondeg}) holds inside any compact set. If we suppose moreover that $ 2 \theta > \sup_x \sigma^2 (x) ,$ then (\ref{lyapunovcondition-5}) is satisfied outside a sufficiently large compact $\tilde K$:   
\begin{eqnarray*}
 2 x^{\top} b(s,x) + tr (a(x))& =&  2 \theta x ( S(s)-x ) + \sigma^2 (x) c_0 + \sigma^2 (x) x^2 \left( 1 - \frac{c_1}{x} \right)^2 \\
  &\le & - x^2 \, [  2 \theta - \sup_x \sigma^2 (x) ] \;+\; o (x^2 ) \quad\mbox{  as  } |x| \to \infty    
\end{eqnarray*} 
by boundedness of $S(s)$. Thus for sufficiently large compact and convex sets, all assumptions of theorem 1.1 hold. Thus we have positive Harris recurrence of the grid chain  {\bf (H)} together with (\ref{nummelincondition-1}). 
 \halmos\\

{\bf Acknowledgments.}  The authors would like to thank Vlad Bally very warmly for a lot of discussions on lower bounds of densities. We are in debt to  Michael Diether who pointed out an error in an earlier version of the manuscript. Ce travail a b\'en\'efici\'e d'une aide de l'Agence Nationale de la Recherche portant la r\'ef\'erence ANR-08-BLAN-0220-01. \\


\hfill {\bf 06.\ 03.\ 2011}
\vskip0,5cm

Reinhard H\"opfner\\
Institut f\"ur Mathematik, Universit\"at Mainz, D--55099 MAINZ, Germany\\
email: {\tt hoepfner@mathematik.uni-mainz.de}
 
\bigskip

Eva L\"ocherbach\\
CNRS UMR 8088, D\'epartement de Math\'ematiques, Universit\'e de Cergy-Pontoise, F--95000 CERGY-PONTOISE,  France\\
email: {\tt eva.loecherbach@u-cergy.fr}

\end{document}